\documentclass[11pt]{article}

\usepackage{latexsym}
\usepackage{amsfonts}
\usepackage[all]{xy}

\newcommand{\diver}{{{\rm{div}\,}}}
\newcommand{\Diver}{{{\rm{Div}\,}}}
\newcommand{\grad}{{{\rm{grad}\,}}}

\newcommand{\Hess}{{{\rm Hess}\,}}
\newcommand{\noin}{\noindent}

\newtheorem{theorem}{\bf Theorem}[section]
\newtheorem{lemma}[theorem]{\bf Lemma}
\newtheorem{corollary}[theorem]{\bf Corollary}
\newtheorem{proposition}[theorem]{\bf Proposition}
\newtheorem{remark}[theorem]{\bf Remark}
\newtheorem{definition}[theorem]{\bf Definition}

\begin{document}

\title{An Extension of Barta's Theorem \\  and Geometric Applications}

\author{G. Pacelli Bessa  \and
 J. F\'{a}bio Montenegro\thanks{\small Both author were partially supported by CNPq Grant.}}
\date{\today}
\maketitle
\begin{abstract}We prove an extension of a theorem of
Barta then
 we make few geometric applications.  We extend
Cheng's
 lower eigenvalue estimates of normal geodesic  balls. We generalize Cheng-Li-Yau
eigenvalue estimates of minimal submanifolds of the space forms.
We prove an stability theorem for minimal hypersurfaces of the
Euclidean space, giving a converse statement of a result of
Schoen. Finally we prove a generalization of a result of
Kazdan-Kramer about existence of solutions of certain quasi-linear
elliptic equations.

\vspace{.1cm}

\vspace{.1cm}

\end{abstract}
\section{Introduction} \hspace{.5cm}The fundamental tone
$\lambda^{\ast}(M)$ of  a smooth Riemannian manifold $M$ is
defined by $$ \lambda^{\ast}(M)=\inf\{\frac{\smallint_{M}\vert
\grad f \vert^{2}}{\smallint_{M} f^{2}};\, f\in
 {L^{2}_{1,0}(M )\setminus\{0\}}\}$$
\noin where $L^{2}_{1,\,0}(M )$ is the completion  of
$C^{\infty}_{0}(M )$ with respect to the norm $ \Vert\varphi
\Vert^2=\int_{M}\varphi^{2}+\int_{M} \vert\nabla \varphi\vert^2$.
 When $M$ is an open   Riemannian manifold,  the
fundamental tone $\lambda^{\ast}(M)$ coincides with the greatest
lower bound  $\inf \Sigma$ of the spectrum $\Sigma \subset
[0,\infty)$ of the unique self-adjoint extension of the Laplacian
$\triangle$ acting on $C_{0}^{\infty}(M)$ also denoted by
$\triangle$. When $M$ is compact with piecewise smooth boundary
$\partial M$ (possibly empty) then $\lambda^{\ast}(M)$ is the
first eigenvalue $\lambda_{1}(M)$ of $M$ (Dirichlet boundary data
if $\partial M\neq \emptyset$). A natural question is what bounds
can one give for $\lambda^{\ast}(M)$ in terms of Riemannian
invariants? Or
 if $M$ is an open manifold, for what geometries does $M$   have $\lambda^{\ast}(M)>0$? See
\cite{kn:berger-gauduchon-mazet}, \cite{kn:chavel}, \cite{kn:G} or
\cite{kn:schoen-yau}. A simple method for giving bounds the first
Dirichlet eigenvalue  $\lambda_{1}(M)$ of a compact smooth
Riemannian manifold $M$ with piecewise smooth
boundary\footnote{\label{foot1}Piecewise smooth boundary here
means that there is a closed set $Q\subset
\partial M$ of $(n-1)$-Hausdorff measure zero such that for each point $q \in
\partial M\setminus Q$ there is a neighborhood of $q$ in $\partial M$ that is a   graph of a smooth function over the
tangent space $T_{q}\partial M$,
 see Whitney \cite{kn:whitney} pages 99-100.
} $\partial M\neq \emptyset$ is Barta's Theorem. All Riemannian
manifolds in this paper are smooth and connected.
\begin{theorem}[Barta, \cite{kn:barta}]\label{barta}Let $M$ be a
bounded
    Riemannian manifold with piecewise smooth non-empty boundary
    $\partial M$ and $f\in C^{2}(M)\cup C^{0}(\overline{M})$ with
    $f\vert M>0$ and $f\vert \partial M=0$ and  $\lambda_{1}(M)$
    be the first Dirichlet eigenvalue of $M$. Then
\begin{equation}\label{eqBarta1}\sup_{M} ( -\triangle f/f) \geq
\lambda_{1}(M ) \geq \inf_{M} (-\triangle f/f).
\end{equation}With equality in (\ref{eqBarta1}) if and only if   $f$
is a first eigenfunction of $M$.
\end{theorem}

\begin{remark}\label{remarkBarta}To obtain the lower bound for $\lambda_{1}(M )$
 we may suppose that  $f\vert\partial M \geq 0$.
\end{remark} Cheng
 applied Barta's Theorem  in a beautiful result
known as Cheng's  eigenvalue comparison theorem.
\begin{theorem}[Cheng, \cite{kn:cheng1}] Let $N$ be a Riemannian
$n$-manifold and $B_{N}(p,r)$ be a geodesic ball centered at $p$
with radius $r<inj (p)$.  Let $c$ be the least upper bound for all
sectional curvatures at $B_{N}(p,r)$ and let $\mathbb{N}^{n}(c)$
be the
 simply connected $n$-space form of constant sectional
curvature $c$. Then
\begin{equation}\label{eqThmcheng1}\lambda_{1}(B_{N}(p, r))\geq
\lambda_{1}(B_{\mathbb{N}^{n}(c)}(r)).
\end{equation}
\end{theorem}In particular, when $c=-1$ and $inj
(p)=\infty$,  Cheng's inequality (\ref{eqThmcheng1}) becomes
$\lambda^{\ast}(N)\geq \lambda^{\ast}(\mathbb{H}^{n}(-1))$ which
is McKean's inequality, see \cite{kn:mckean}. Our main result
 gives lower bound for the fundamental tone $\lambda^{\ast}(M)$ of an arbitrary smooth Riemannian
manifolds $M$ in terms of the divergence of certain vector fields
regardless  the smoothness degree of its boundary.

\begin{definition}Let $M$ be a Riemannian
manifold and a vector field $X\in L^{1}_{loc}(M)$ $($meaning that
$\vert X\vert \in L_{loc}^{1}(M))$. A function $g\in
L_{loc}^{1}(M)$ is a weak divergence of $X$ if
$$\int_{M} \phi \,g =-\int_{M}\langle \grad
\phi,\,X\rangle,\,\, \forall \,\phi \in
C^{\infty}_{0}(M).$$\label{defDiv1}\end{definition}
 It is clear that there exists at
most one weak divergence $g \in L^{1}_{loc}(M)$ for a given $X\in
L^{1}_{loc}(M)$  and we may write $g=\Diver X$.  For $C^{1}$
vector fields $X$ the classical divergence $\diver X$ and the weak
divergence $\Diver X$ coincide.
\begin{definition}\label{defDiv2}Let ${\cal W}^{1,1}(M)$ denote the Sobolev space  of all vector
fields $ X \in L^{1}_{loc}(M)$ possessing weak divergence $ \Diver
X $.
\end{definition}
\begin{remark}\label{remarkW11} If $X \in {\cal
W}^{1,1}(M)$ and $f\in C^{1}(M)$ then $fX\in {\cal W}^{1,1}(M)$
with $\Diver (f X) = \langle \grad f,\,X\rangle + f \,\Diver X$.
In particular for  $f \in C^{\infty}_{0}(M)$ we have that
\begin{equation}\label{eqDivergence1}\int_{M}\Diver (fX)=\int_{M} \langle \grad
f,\,X\rangle -  \langle \grad
f,\,X\rangle=0.\end{equation}Conversely,  if $fX\in {\cal
W}^{1,1}(M)$ for all $f\in C_{0}^{\infty}(M)$ then $X\in {\cal
W}^{1,1}(M)$.
\end{remark}

Our main result is the following theorem.
\begin{theorem}\label{teoremaPrincipal}Let
$M$ be  a Riemannian manifold. Then
\begin{equation}\label{eqThmP1}
 \lambda^{\ast}(M)\geq
 \sup_{{\cal W}^{1,1} }\{\inf_{ M} (\Diver X-
\vert X \vert^{2})\}. \end{equation} If $M$ is a compact
Riemannian manifold with smooth boundary then
\begin{equation}\label{eqThmP2}\lambda_{1}(M)=
 \sup_{{\cal W}^{1,1} }\{\inf_{ M} (\Diver X-
\vert X \vert^{2})\} .\end{equation}
\end{theorem}
Our first geometric application of Theorem
(\ref{teoremaPrincipal}) is an extension of Cheng's lower
eigenvalue estimates. We show that inequality (\ref{eqThmcheng1})
is valid for arbitrary geodesic balls $B_{N}(p,r)$ provided the
$(n-1)$-Hausdorff measure ${\cal H}^{n-1}({\rm Cut} (p)\cap
B_{N}(p,r))=0$, where ${\rm Cut}(p)$ is the cut locus of $p$.
Moreover, we show that equality of the eigenvalues  occurs if and
only if $B_{N}(p,r)$ and
 $B_{\mathbb{N}^n(c)}(r)$ are isometric.

\begin{theorem} \label{Cheng} Let $N$ be a  Riemannian $n$-manifold  with radial sectional
curvature $K(x)(\partial t, v)\leq c$, $x\in B_{N}(p,r)\setminus
{\rm Cut}(p)$,  $v\perp \partial t$ and  $\vert v\vert \leq 1$.
Let $\mathbb{N}^{n}(c)$ be the
 simply connected $n$-space form of constant sectional
curvature $c$ and suppose that ${\cal H}^{n-1}({\rm Cut} (p)\cap
B_{N}(p,r))=0$. Then
\begin{equation}\label{eqCheng}\lambda^{\ast}(B_{N}(p,r))\geq
\lambda_{1}(B_{\mathbb{N}^{n}(c)}(r)).
\end{equation}Equality in (\ref{eqCheng}) holds iff
$B_{N}(p,r)$ and $B_{\mathbb{N}^{n}(c)}(r)$ are isometric.
\end{theorem} Our second geometric application (Theorem \ref{submanifold}) is a generalization
of the following
   Cheng-Li-Yau's result proved in
\cite{kn:cheng-li-yau}.
\begin{theorem}[Cheng-Li-Yau, \cite{kn:cheng-li-yau}] Let
$M^{m}\subset \mathbb{N}^{n}(c)$ be an immersed $m$-dimensional
minimal
  submanifold   of the  $n$-dimensional space form of constant sectional curvature $c\in\{-1,0,1\}$
   and $D\subset M^{m}$ a $C^{\,2}$ compact domain. Let $a=\inf_{p\in D}\sup_{z\in D}
   dist_{\mathbb{N}^{n}(c)}(p,z)>0$ be the outer radius of $D$. If $c=1$ suppose that $a\leq
   \pi/2$. Then
   \begin{equation}\label{eqCLY1}\lambda_{1}(D)\geq
   \lambda_{1}(B_{\mathbb{N}^{m}(c)}(a))
   \end{equation}
   Equality in (\ref{eqCLY1}) holds iff $M$ is totally geodesic in
   $\mathbb{N}^{n}(c)$ and $D=B_{\mathbb{N}^{m}(c)}(a)$.
   \end{theorem}
\begin{theorem}\label{submanifold}Let $N$ be a  Riemannian $n$-manifold  with radial sectional
curvature $K(x)(\partial t, v)\leq c$, for all $x\in
B_{N}(p,r)\setminus {\rm Cut}(p)$, and all $v\perp \partial t $
with $\vert v\vert \leq 1$. Let $M\subset N$ be an $m$-dimensional
minimal submanifold and $\Omega \subset M\cap B_{N}(p,r)$ be a
connected component. Suppose that the $(m-1)$-Hausdorff measure
${\cal H}^{m-1}(\Omega \cap {\rm Cut}_{N}(p))=0$. If $c>0 $,
suppose in addition that $r< \pi /2\sqrt c$. Then
\begin{equation}\label{eqSubm1}\lambda^{\ast}(\Omega) \geq
\lambda_{1}(B_{\mathbb{N}^{m}(c)}(r)),
\end{equation}where $B_{\mathbb{N}^{m}(c)}(r)$  is the geodesic
ball with radius $r$ in the simply connected space form
$\mathbb{N}^{n}(c)$ of constant sectional curvature $c$. If
$\Omega $ is bounded then equality in (\ref{eqSubm1}) holds iff
$\Omega=B _{\mathbb{N}^{m}(c)}(r)$ and $M=\mathbb{N}^{m}(c)$.
\end{theorem}
 After Nadirashvili 's  bounded
  minimal surfaces in $\mathbb{R}^{3}$, (see \cite{kn:Nadirashvili}),  Yau in \cite{kn:yau} asked if  the spectrum of
a Nadirashvili minimal surface  was discrete. A more basic
question is if the lower bound of the spectrum of a Nadirashvili
minimal surface is positive. The following corollary
 shows that this is the case.
 \begin{corollary}\label{nadirashvili}Let $M\subset B_{\mathbb{R}^{3}}(r)$ be a
 complete bounded minimal surface in $\mathbb{R}^{3}$. Then
 $$\lambda^{\ast}(M)\geq
 \lambda_{1}(\mathbb{D}(r))=c/r^{2}.
$$ \noin
 Where $c>0$ is an absolute constant.
 \end{corollary}
 Let
$M\subset\mathbb{R}^{3}$ be  a minimal  surface with second
fundamental form $A$ and $B_{M}(p,r)$  be a stable geodesic ball
  with radius $r$. Schoen,  \cite{kn:schoen} showed
that  $\Vert A \Vert^{2}(p) \leq c/r^{2}$ for an absolute constant
$c>0$. We have a converse of Schoen's result for minimal
hypersurfaces of the Euclidean space.
\begin{corollary}\label{stability}Let $M\subset \mathbb{R}^{n+1}$ be
  a minimal hypersurface  with second fundamental form
$A$ and $B_{M}(r)$  be a  geodesic ball   with radius $r$. If
$$\sup_{B_{M}(r)} \Vert A \Vert^{2} \leq
\lambda_{1}(B_{\mathbb{R}^{n}}(0,r))=c(n)/r^2, $$ \noin then
$B_{M}(r)$ is stable. Here $B_{\mathbb{R}^{n}}(0,r)$ is a geodesic
ball of radius $r$ in the Euclidean space $\mathbb{R}^{n}$ and
$c(n)>0$ is a constant depending on $n$.
\end{corollary}
Finally we apply Barta's Theorem and Theorem
(\ref{teoremaPrincipal}) to the theory of quasi-linear elliptic
equations.
\begin{theorem}\label{Elliptic1} Let $M$ be a bounded Riemannian manifold with smooth
boundary and
 $F\in C^{0}(\overline{M})$.
Consider this problem,
\begin{equation}\label{eqElliptic1}\left\{\begin{array}{rccrl} \triangle u - \vert \grad
u\vert^{2} &=& F & in&
M\\
u&=& +\infty & on& \partial M .
\end{array}\right.
\end{equation}If (\ref{eqElliptic1})  has a smooth
solution then $\inf_{M} F\leq \lambda_{1}(M) \leq \sup_{M} F $. If
either $\inf_{M} F=\lambda_{1}(M)$ or $\lambda_{1}(M)= \sup_{M} F$
then $F=\lambda_{1}(M)$. On the other hand if
 $F=\lambda$ is a constant  then the problem (\ref{eqElliptic1}) has solution if
and only if $\lambda=\lambda_{1}(M)$.
\end{theorem}
 Now if  we
allow continuous boundary data on problem (\ref{eqElliptic1}) we
have the following generalization of a result of Kazdan-Kramer
\cite{kn:kazdan-kramer}.
\begin{theorem}\label{Elliptic2} Let $M$ be a bounded Riemannian manifold with smooth
boundary and
 $F\in C^{0}(\overline{M})$ and $\psi \in C^{0}(\partial
M)$. Consider the problem
\begin{equation}\label{eqElliptic2}\left\{\begin{array}{rclll} \triangle u
- \vert \grad u\vert^{2} &=& F & in&
M\\
u&=& \psi & on& \partial M.
\end{array}\right.
\end{equation}then if $\sup_{M}F < \lambda_{1}(M)$ then (\ref{eqElliptic2}) has
solution. Moreover if (\ref{eqElliptic2}) has
 solution then $\inf_{M}F< \lambda_{1}(M)$.
\end{theorem}
\vspace{.1cm}

\begin{remark}\label{remarkkazdan}
~
\begin{itemize}
\item[1] If we set $f=e^{-u}$ then   (\ref{eqElliptic1}) becomes
\begin{equation}\label{eqremark2}\left\{\begin{array}{rclll}\triangle \,f + F\,f&= &0 & in&
M\\
 f&=& 0 & on& \partial M .
\end{array}\right.
\end{equation}
Kazdan-Warner in \cite{kn:kazdan-warner}  studied this problem
(\ref{eqremark2}) and they  showed that if $F\leq \lambda_{1}(M)$,
then (\ref{eqremark2}) has solution, with $\sup F=\lambda_{1}(M)$.
Thus Theorem (\ref{Elliptic1}) is a complementary result to
Kazdan-Warner result.
 \item[2]If we impose Dirichlet boundary data
($u=0$ on $\partial M$) on problem (\ref{eqElliptic1}) then  there
is a solution if $F\leq \lambda_{1}(M)$ with strict inequality in
a positive measure subset of $M$. This was proved by Kazdan-Kramer
in \cite{kn:kazdan-kramer}.
\end{itemize}
\end{remark}
 \noin
{\bf Acknowledgement:}{\em The beginning of this work was greatly
motivated by a paper of Kazdan \& Kramer \cite{kn:kazdan-kramer}
for which we thank Professor Djairo de Figueiredo for bringing it
to our attention. We also should mention a paper of Cheung-Leung
\cite{kn:cheung-leung} which contains the germ of Theorem
(\ref{teoremaPrincipal}). We are grateful to the referee for
pointing out a gap in the proof of Theorem (\ref{Cheng}).}

\section{An extension of Barta's theorem}
\noin We can consider Theorem (\ref{teoremaPrincipal}) as an
extension of Barta's theorem for if $M$ is a bounded Riemannian
manifold with piecewise smooth boundary $\partial M\neq \emptyset$
and $f\in C^{\,2}(M)\cup C^{\,0}(\overline{M})$ is a positive
function on $M$ then  for the vector field $X=-\grad \log f$  we
obtain that $\diver X-\vert X \vert^{2}=-\triangle f/f$.

\vspace{.5cm}
 \noin {\bf
Theorem \,\ref{teoremaPrincipal}} {\em\, Let $M$ be a Riemannian
manifold. Then
\begin{equation}\label{eqThmP3}
 \lambda^{\ast}(M)\geq
 \sup_{{\cal W}^{1,1}(M ) }\{\inf_{ M} (\Diver X-
\vert X \vert^{2})\}.
\end{equation}
If $M$ is compact with smooth non-empty boundary then
\begin{equation}\label{eqThmP4}\lambda_{1}(M)=
 \sup_{{\cal W}^{1,1}(M ) }\{\inf_{ M} (\Diver X-
\vert X \vert^{2})\}. \end{equation} } \noin {\bf Proof:}  Let
$X\in {\cal W}^{1,1}(M) $ and $f\in C_{0}^{\infty}(M)$. As
observed in the Remark (\ref{remarkW11}) we have that
$\smallint_{M}\Diver (f^{2} X)=0$. On the other hand we
 have that

\begin{eqnarray}0=\int_{M}\Diver(f^{2}X) & = & \int_{M} \langle grad\,f^{2}, X\rangle
   + \int_{M} f^{2}\,\Diver (X) \nonumber \\
            & & \nonumber\\
                            &\geq & -\int_{M}\vert grad\,f^{2}\vert\cdot \vert X \vert+
                             \, \int_{M}f^{2}\,\Diver\,X
                              \nonumber\\
                             & &  \nonumber \\
                            &= & -\int_{M}2 \cdot \vert f\vert\ \cdot   \,\vert X \vert
\cdot \vert grad\,f\vert + \,
\int_{M}f^{2}\,\Diver\,X  \nonumber\\
 & &  \nonumber \\
 & \geq & -\int_{M} f^{2}\cdot \vert X \vert^{2} - \int_{M}\vert \grad f \vert^{2}
 +\int_{M}f^{2}\,\Diver\,X \nonumber\\
  &&  \nonumber\\
   & = & \int_{M}(\Diver\,X -\vert X \vert^{2})\cdot f^{2} - \int_{M}\vert \grad f
   \vert^{2} \nonumber\\
   && \nonumber \\
    &\geq & \inf_{ M}(\Diver\,X -\vert X
    \vert^{2})\int_{M}f^{2} - \int_{M}\vert \grad f
   \vert^{2}. \nonumber
\end{eqnarray}
Then $$ \int_{M}\vert \grad f \vert^{2} \geq   \inf_{ M}(\Diver\,X
-\vert X
    \vert^{2})\int_{M}f^{2},
$$ and thus
$$
\int_{M}\vert \grad f
   \vert^{2}  \geq  \sup_{{\cal W}^{1,1} }\inf_{ M}(\Diver\,X -\vert X
    \vert^{2})\int_{M}f^{2}.
$$\noin
Therefore $$\lambda^{\ast}(M) \geq  \sup_{{\cal W}^{1,1} }\inf_{
M}(\Diver\,X -\vert X
    \vert^{2}).$$ \noin This proves (\ref{eqThmP3}).
Suppose that $M$ is compact with smooth non-empty boundary  and
let $v$ be its first eigenfunction. If we set $X_{0}=-\grad (\log
v)$ then we have that $\diver X_{0} -\vert X_{0}
    \vert^{2}= -\triangle v/v=\lambda_{1}(M)$ and (\ref{eqThmP4}) is
    proven.

\begin{remark}\label{remarkBartaGen}This same proof above shows that
$$\lambda^{\ast}(M) \geq  \sup_{{\cal W}^{1,1} }\inf_{ M\setminus
G}(\Diver\,X -\vert X
    \vert^{2}),$$ where $G$ has zero Lebesgue measure.
\end{remark}
\section{Geometric Applications}
In the geometric applications  of the Theorem
(\ref{teoremaPrincipal}) we need to know when a given vector field
$X$ is in the Sobolev space ${\cal W}^{1,1}(M)$. The following
lemma give sufficient conditions.
\begin{lemma}\label{divergence}
Let $\Omega\subset M$ be a bounded domain with  piecewise smooth
boundary $\partial \Omega$  in a smooth  Riemannian manifold $M$.
Let $G\subset \Omega$ be a closed subset with $(n-1)$-Hausdorff
measure ${\cal H}^{n-1}(G)=0$. Let $X$ be a vector field  of class
$ C^{1}(\Omega\setminus G)\cap L^{\infty}(\Omega)$ such that
$\diver (X)\in L^{1}(\Omega)$. Then $X \in {\cal
W}^{1,1}(\Omega)$.
\end{lemma}
\noin {\bf Proof:} In fact we are going to show that
\begin{equation}\label{eqDivergence2}\int_{\Omega}\diver (X)
= \int_{\partial \Omega}\langle X,\nu\rangle,
\end{equation}where $\nu$ is the outward unit normal vector field
on $\partial \Omega \setminus Q$ and $Q\subset \partial \Omega$ is
a closed subset with ${\cal H}^{n-1}(Q)=0$, see the footnote
$(\ref{foot1})$ on page $2$. The equation (\ref{eqDivergence2})
implies that $X\in {\cal W}^{1,1}(\Omega)$ since for all $\phi \in
C^{\infty}_{0}$ we have that
$$\int_{\Omega}\diver (\phi\,X)=\int_{\Omega}\phi\,\diver (X)+\int_{\Omega}\langle \grad \phi ,X \rangle
$$ and $\int_{\Omega}\diver (\phi\,X)=0$ by
(\ref{eqDivergence2}). Suppose first that $M=\mathbb{R}^{n}$. We
may assume   that $G$ is connected otherwise we work with each
connected component. By Whitney's Theorem there exists a
 smooth function (we may suppose to be non-negative)
$f\hspace{-.7mm}:\hspace{-.7mm}\mathbb{R}^{n}\hspace{-.7mm}
\rightarrow \mathbb{R}$ such that $G\hspace{-.8mm}=f^{-1}(0)$. By
Sard's theorem we can pick a sequence of regular values
$\epsilon_{i}\rightarrow 0$ such that $f^{-1}(\epsilon_{i})$ is a
smooth $(n-1)$-dimensional submanifold $\partial N_{i}$ bounding a
connected set $N_{i}=f^{-1}([0,\epsilon_{i}])$ containing $G$,
moreover, $N_{i}\subset \Omega$ for $\epsilon_{i}$ sufficiently
small. Set $\Omega_{i}=\Omega \setminus N_{i}$ and let $\chi_{i}$
its characteristic function. Then $
 \chi_{i}\cdot \diver X\hspace{-.1cm} \to \diver X$ a.e. in $\Omega$ and $\chi_{i}\cdot \diver X \leq\diver X\in
L^{1}(\Omega)$. By the Lebesgue  Convergence Theorem
$\int_{\Omega_{i}}\diver X=\int_{\Omega}\chi_{i}\cdot \diver X \to
\int_{\Omega} \diver X $. On the other hand applying the
divergence theorem to $X$ on $\Omega_{i}$ we obtain
$$\begin{array}{lcl}
\int_{\Omega_{i}}\diver X & = & \int_{\partial \Omega_{i}}\langle X, \nu \rangle\\
&& \\
                     &=& \int_{\partial \Omega}\langle X, \nu \rangle -\int_{\partial N_{i}}\langle X, \nu_{i}
                     \rangle ,
\end{array}
$$ $\nu_{i}$ is the outward (pointing toward $G$) unit vector field normal to $\partial N_{i} $.
 But $\left \vert \int_{\partial N_{i}}\langle X, \nu_{i}
\rangle\right\vert \leq
 Vol_{n-1}(\partial N_{i})  \Vert X\Vert_{\infty}=Vol_{n-1}(f^{-1}(\epsilon_{i}))  \Vert
 X\Vert_{\infty}$ and $Vol_{n-1}(f^{-1}(\epsilon_{i}))\to {\cal
 H}^{n-1}(G)=0$ as we will show later.
Passing to the limit we have
 $$\begin{array}{lcl}\int_{\Omega} \diver X& = & \lim_{i\to \infty}\int_{\Omega_{i}}\diver X\\ && \\& =&
\int_{\partial \Omega}\langle X, \nu \rangle
-  \lim_{i\to \infty}\int_{\partial N_{i}}\langle X, \nu_{i} \rangle\\
&& \\
&=& \int_{\partial \Omega}\langle X, \nu \rangle\end{array}$$ To
show $Vol_{n-1}(f^{-1}(\epsilon_{i}))\to {\cal
 H}^{n-1}(G)=0$
recall the $(n-1)$-dimensional spherical measure of $A\subset
\mathbb{R}^{n}$ is  defined by ${\cal
S}^{n-1}(A)=\sup_{\delta>0}{\cal S}^{n-1}_{\delta}(A)$
$=\lim_{\delta\to 0}{\cal S}^{n-1}_{\delta}(A)$, where ${\cal
S}^{n-1}_{\delta}(A)=\inf\sum {\rm diam}
 (B_{j})^{n-1}$ the infimum taken over all coverings of $A$
  by balls $B_{j}$ with diameter  ${\rm diam} (B_{j})\leq \delta$.
 The
 $(n-1)$-spherical measure is related to $(n-1)$-Hausdorff measure
 by  ${\cal
H}^{n-1}(A)\leq {\cal S}^{n-1}(A)\leq 2^{n-1}{\cal H}^{n-1}(A)$,
$\forall \,A\subset \mathbb{R}^{n}$, see \cite{kn:matilla} for
more details. Since ${\cal S}^{n-1}(G)=0$ we have that ${\cal
S}^{n-1}_{\delta}(G)=0$ for all $\delta>0$. Then for $\delta>0$
and all $k>0 $ there exists a (finite) covering of $G$, ($G$ is
compact),
 by closed balls  $\{B_{kj}\}$ of diameter ${\rm
diam}(B_{kj})\leq\delta$ such that $\sum_{j} {\rm diam}
(B_{kj})^{n-1} \leq 1/k$.  For $\epsilon_{i}$ sufficiently small,
say $\epsilon_{i}\leq\epsilon_{i_{0}}$, the balls $B_{kj}$ are
also a covering for the submanifold $f^{-1}(\epsilon_{i})$. This
means that ${\cal S}^{n-1}_{\delta}(f^{-1}(\epsilon_{i}))\leq
1/k$. Taking a sequence $\delta_{l}\to 0$ we can find a sequence
of pairs $(k_{l},\epsilon_{l})$,
 $k_{l}\rightarrow
\infty$, $\epsilon_{l}\to 0$ such that
 ${\cal S}_{\delta_{l}}^{n-1}(f^{-1}(\epsilon_{i}))\leq 1/k_{l}$ for all $\epsilon_{i}\leq
 \epsilon_{l}$. Thus we have that
$\lim_{l\to \infty}{\cal S}^{n-1}_{\delta_{l}}
(f^{-1}(\epsilon_{l}))=0$ and since ${\cal S}^{n-1}
(f^{-1}(\epsilon_{l}))\leq {\cal S}^{n-1}_{\delta_{l}}
(f^{-1}(\epsilon_{l}))$ we have that $\lim_{l\to \infty}{\cal
S}^{n-1} (f^{-1}(\epsilon_{l}))=0$.

\noin  On the other
 ${\cal S}^{n-1}
(f^{-1}(\epsilon_{l}))\geq {\cal
H}^{n-1}(f^{-1}(\epsilon_{l}))=c(n)\cdot
vol_{n-1}(f^{-1}(\epsilon_{l}))$. This shows that
$vol_{n-1}(f^{-1}(\epsilon_{i}))\to 0$. This completes the proof
of Lemma (\ref{divergence}) for $M=\mathbb{R}^{n}$. The general
case is done similarly  using partition of unit. \vspace{.5cm}

\subsection{Geodesic coordinates}
  Let $M$ be a  Riemannian manifold and a point $p\in M$.  For
  each vector
  $\xi \in T_{p}M$, let $\gamma_{\xi}$ the unique geodesic satisfying
  $\gamma_{\xi}(0)=p$,
  $\gamma_{\xi}'(0)=\xi$ and $d(\xi)=\sup\{t>0:
  {\rm dist}_{M}(p,\gamma_{\xi}(t))=t\}$. Consider the largest open subset
  ${\cal D}_{p}=\{t\,\xi\in T_{p}M: 0\le t<d(\xi), \,\vert
   \xi\vert=1\}$ of
  $T_{p}M$    such that for
  any $\xi \in {\cal D}_{p}$ the  geodesic
  $\gamma_{\xi}(t)=\exp_{p}(t\,\xi)$ minimizes the distance from
  $p$ to $\gamma_{\xi}(t)$ for all $t\in [0,d(\xi)]$. The cut locus of $p$ is given by ${\rm Cut }(p)=\{ \exp_{p}(d(\xi)\,\xi),\, \xi \in T_{p}M,\, \vert
  \xi \vert =1\}$ and  $M=\exp_{p}({\cal D}_{p})\cup {\rm Cut
  }(p)$. The exponential map $\exp_{p}:{\cal D}_{p}\to \exp_{p}({\cal
  D}_{p})$ is a diffeomorphism and is called the geodesic
  coordinates of $M\setminus {\rm Cut}(p)$.

\noin Fix a vector $\xi \in T_{p}M$, $\vert \xi \vert =1$ and
   denote by $\xi^{\perp}$  the orthogonal
complement of $\{\mathbb{R}\xi\}$ in $T_{p}M$ and let
$\tau_{t}:T_{p}M\to T_{\exp_{p}(t\,\xi)}M$ be the parallel
translation along $\gamma_{\xi}$.  Define the path of linear
transformations $$ {\cal A}(t,\xi):\xi^{\perp}\to\xi^{\perp}$$ by
$${\cal A}(t,\xi)\eta=(\tau_{t})^{-1}Y(t)$$ where $Y(t)$ is the
Jacobi field along $\gamma_{\xi}$ determined by the initial data
$Y(0)=0$, $(\nabla_{\gamma_{\xi}'}Y)(0)=\eta$. Define the map
$${\cal R}(t):\xi^{\perp}\to \xi^{\perp} $$  by $${\cal
R}(t)\eta=(\tau_{t})^{-1}\, {\rm
R}(\gamma_{\xi}'(t),\tau_{t}\,\,\eta)\gamma_{\xi}'(t),$$  where
${\rm R}$ is the Riemann curvature tensor of $M$. It turns out
that the map ${\cal R}(t)$ is a self adjoint map and the path of
linear transformations ${\cal A}(t,\xi)$ satisfies the Jacobi
equation ${\cal A}''+{\cal R}{\cal A}=0$ with initial conditions
${\cal A}(0,\xi)=0$, ${\cal A}'(0,\xi)=I$. On the set
$\exp_{p}({\cal D}_{p})$ the Riemannian metric of $M$ can be
expressed by
\begin{equation}\label{eqmetricGeoCoord}ds^{2}(\exp_{p}(t\,\xi))=dt^{2}+\vert {\cal
A}(t,\xi)d\xi\vert^{2}. \end{equation} Setting
$\sqrt{g(t,\xi)}=\det{\cal A}(t,\xi)$ we have by Rauch comparison
theorem this following comparison theorem due to R. Bishop
\cite{kn:bishop-crittenden}, see also \cite{kn:chavel}.
\begin{theorem}[Bishop, \cite{kn:bishop-crittenden}]\label{thmBishop}If
the radial sectional curvatures along $\gamma_{\xi}$ satisfies
$\langle {\rm R}(\gamma_{\xi}'(t),v)\gamma_{\xi}'(t),v\rangle \leq
c\vert v\vert^{2}$, $\forall t\in (0,r)$ and if $S_{c}(t)$ does
not vanishes on $(0,r)$ then
\begin{eqnarray}\label{eqRauch}
[\sqrt{g(t,\xi)}/S_{c}^{n-1}(t)]'&\geq &0\,\,on \,\,(0,r)\nonumber \\
&& \\
 \sqrt{g(t,\xi)}-S_{c}^{n-1}(t)&\geq&0\,\,on \,\,(0,r]\nonumber \end{eqnarray}
Moreover  equality occurs in one of these two inequalities
(\ref{eqRauch}) at a point $t_{0}\in (0,r)$ iff
 ${\cal R}=c\cdot I$ and ${\cal A}=S_{c}\cdot I $ on all
$[0,t_{0}]$.\end{theorem}
 Where $C_{c}(t)=S_{c}'(t)$ and $S_{c}(t)$ is
given by
 \begin{equation}\label{eqSc}S_{c}(t)=\left\{\begin{array}{clll}\displaystyle\frac{1}{\sqrt c}\sin
 (\sqrt c \,t)
 , & if& c>0 &\\
 t & if & c=0\\
\displaystyle\frac{1}{\sqrt-c} \sinh ( \sqrt-c\, t) & if & c<0
 \end{array}\right.
 \end{equation}

 \subsection{Proof of Theorem \ref{Cheng}}

 \noin {\bf Theorem \ref{Cheng}} {\em   Let $N$ be a
Riemannian $n$-manifold  with radial sectional curvature
$K(x)(\partial t, v)\leq c$,  $x\in B_{N}(p,r)\setminus {\rm
Cut}(p)$ and $v\in T_{x}N\cap (\partial t)^{\perp}$ with $\vert
v\vert \leq 1$. Let $\mathbb{N}^{n}(c)$ be the
 simply connected $n$-space form of constant sectional
curvature $c$ and suppose that ${\cal H}^{n-1}({\rm Cut} (p)\cap
B_{N}(p,r))=0$. Then
\begin{equation}\label{eqCheng2}\lambda^{\ast}(B_{N}(p,r))\geq
\lambda_{1}(B_{\mathbb{N}^{n}(c)}(r)).
\end{equation}Equality in (\ref{eqCheng2}) holds iff
$B_{N}(p,r)$ and $B_{\mathbb{N}^{n}(c)}(r)$ are isometric.}

\vspace{.2cm} \noin {\bf Proof:} \noin Observe that if $c>0$ and
$r>\pi/\sqrt{c}$ there is nothing to prove. Because in this case
$\mathbb{N}^{n}(c)= \mathbb{S}^{n}(c)=B_{\mathbb{N}^{n}(c)}(r)$
and $\lambda_{1}(B_{\mathbb{N}^{n}(c)}(r))=0$. Hence, we may
assume that $r<\pi/\sqrt{c}\,$ if $c>0$.  Let $v$ be a positive
first eigenfunction of $B_{\mathbb{N}^{n}(c)}(r)$. It is well
known that $v$ is a radial function satisfying the differential
equation
\begin{equation}v''(t) +
(n-1)\frac{C_{c}(t)}{S_{c}(t)}v'(t)+\lambda_{1}(B_{\mathbb{N}^{n}(c)}(r))v(t)=0,
\,\, \forall \,\,t\in[0,r] \label{eqCheng3}
\end{equation} with $v'(t)\leq 0$ and with
$v'(t)=0$ iff $t=0$.  Define  $u :B_{N}(p,r) \rightarrow
[0,\infty)$ by
$$ u(x)=\left \{ \begin{array}{cll} v(t) & if &
x=\exp_{p}(t\,\xi),\;t\in [0,d(\xi))\cap[0,r]\\
&&\\
0 &if &x=\exp_{p}(d(\xi)\,\xi)\in {\rm Cut} (p)
\end{array}\right.
$$
Set $X(x)=-\grad \log u (x)$ if $x\in
B_{N}(p,r)\setminus(\{p\}\cup {\rm Cut}(p))$ and $X(x)=0$ if $x\in
B_{N}(p,r)\cap (\{p\}\cup {\rm Cut}(p))$. The vector field $X$ is
expressed in geodesic coordinates by
$$X(x)=\left\{\begin{array}{cll}\displaystyle
-\frac{v'(t)}{v(t)}\cdot \partial t & if &x=\exp_{p}(t\,\xi),\;t\in (0,d(\xi))\cap (0,r]\\
& &\\
0&if & x=p,\,\;\;{\rm or}\;\; x=\exp_{p}(d(\xi)\,\xi).
\end{array}\right.
$$ The vector field $X\in {\cal
W}^{1,1}(B_{N}(p,r))$    as we will prove it later. Setting
$B\setminus G=B_{N}(p,r)\setminus (\{p\}\cup {\rm Cut}(p))$ for
simplicity of notation we have by Theorem (\ref{teoremaPrincipal})
and Remark (\ref{remarkBartaGen}) we have that
$$\lambda^{\ast}(B_{N}(p,r)) \geq  \inf_{B\setminus G} \{\Diver X - \vert X
 \vert^{2}\}=\inf_{B\setminus G} \{\diver X - \vert X
 \vert^{2}\}=\inf_{B\setminus G}-\frac{\triangle u}{u}$$

\noin By  (\ref{eqRauch}) and (\ref{eqCheng3}) we have that for
$0<t<d(\xi)$ \begin{eqnarray}-\frac{\triangle
u}{u}(\exp_{p}(t\,\xi))&=&-\displaystyle\frac{1}{v(t)}\left\{v''(t)
+ \frac{\sqrt{g(t,\xi)}'}{\sqrt{g(t,\xi)}} v'(t)\right\}\nonumber
\\
&& \nonumber\\
 &\geq&  -\frac{1}{v(t)}\left\{v''(t) +(n-1)
\frac{S'_{c}(t)}{S_{c}(t)}v'(t) \right\}\label{eqCheng4}
\\ \nonumber \\
& =&
\lambda_{1}(B_{\mathbb{N}^{n}(c)}(r))\nonumber\end{eqnarray}
Thus
$$\lambda^{\ast}(B_{N}(p,r))\geq
\inf_{B\setminus G}[-\displaystyle\frac{\triangle u}{u}]\geq
\lambda_{1}(B_{\mathbb{N}^{n}(c)}(r))$$ This proves
(\ref{eqCheng2}).
 To handle the equality case in (\ref{eqCheng2}) we
observe that the set of non-smooth points of the boundary
$\partial B_{N}(p,r)$ is exactly the intersection   $ {\rm
Cut}(p)\cap\partial B_{N}(p,r)$ but since ${\cal H}^{n-1}(  {\rm
Cut}(p)\cap B_{N}(p,r))=0$ the Hausdorff measure   ${\cal
H}^{n-1}({\rm Cut}(p)\cap
\partial B_{N}(p,r))=0$. Thus $B_{N}(p,r)$ has piecewise smooth
boundary (see \cite{kn:whitney} pages 99-100). Hence there exists
 a positive Dirichlet eigenfunction for
$B_{N}(p,r)$ with eigenvalue
$\lambda_{1}(B_{N}(p,r))=\lambda^{\ast}(B_{N}(p,r))$. The equality
$\lambda^{\ast}(B_{N}(p,r))=
\lambda_{1}(B_{\mathbb{N}^{n}(c)}(r))$  implies  that $u$ is a
 first eigenfunction of $B_{N}(p,r)$, see Barta's Theorem (\ref{barta}). Looking at (\ref{eqCheng4})
 we
 see that the
equality implies that
$$
\frac{\sqrt{g(t,\xi)}'}{\sqrt{g(t,\xi)}}= (n-1)
\frac{S'_{c}(t)}{S_{c}(t)}
$$ on all $[0,r]$. Thus by Bishop's Theorem (\ref{thmBishop}) we have that  ${\cal R}=c\cdot I$ and ${\cal
A}=S_{c}\cdot I $. This is saying that $B_{N}(p,r)$ is isometric
to $ B_{\mathbb{N}^{n}(c)}(r)$.
\vspace{.3cm}

 To finish the proof
of Theorem (\ref{Cheng}) we need to show that the vector field
$X\in {\cal W}^{1,1}(B_{N}(p,r)) $. As observed in Remark
(\ref{remarkW11}) it suffices to show that $f^{2}X\in {\cal
W}^{1,1}(B_{N}(p,r))$ for all $f\in C_{0}^{\infty}(B_{N}(p,r))$.
Since the $(n-1)$-Hausdorff measure ${\cal H}^{n-1}(\{p\}\cup {\rm
Cut}(p))=0$ by the Lemma (\ref{divergence})
 it is sufficient to show that
\begin{itemize}
\item[(i)]$\diver (f^{2}X ) \in L^{1}(B_{N}(p,r))$
\item[(ii)]$f^{2}X\in \ C^{1}\left(B_{N}(p,r)\setminus \{p\}\cup
{\rm Cut}(p)\right)\cap L^{\infty}(B_{N}(p,r))$.
\end{itemize}
The vector field $X$ is clearly smooth on $B_{N}(p,r)\setminus
\{p\}\cup {\rm Cut}(p)$ thus on this set $\diver X = \Diver X$ and
$\diver (f^{2}X)=\langle \grad f^{2},X\rangle + f^{2}\diver X.$
Integrating over  $B_{N}(p,r)\cap {\rm Supp}(f)$ we have
\begin{eqnarray}\label{eqDiv(fX)}\int_{B_{N}(p,r)\cap {\rm Supp}(f)}\vert \diver
(f^{2}X)\vert & \leq & \int_{B_{N}(p,r)\cap {\rm Supp}(f)}\vert
\grad f^{2}\vert \vert X\vert \nonumber \\ && \\ &+ &
\int_{B_{N}(p,r)\cap {\rm Supp}(f)}\vert f^{2}\vert \vert\diver
X\vert\nonumber\end{eqnarray} The first term of (\ref{eqDiv(fX)})
is finite
\begin{eqnarray} \int_{B_{N}(p,r)\cap {\rm Supp}(f)}\vert \grad
f^{2}\vert \vert X\vert & \leq &\sup_{ {\rm Supp}(f)}\{\vert \grad
f^{2}\vert\cdot \vert X\vert\}\cdot vol (B_{N}(p,r))\nonumber \\
&& \nonumber \\ & < & \infty \nonumber
\end{eqnarray} since $\vert X (x) \vert \leq \vert (v'/v)(t)\vert  <\infty
$ for $x=\exp_{p}(t\,\xi)\in B_{N}(p,r)\cap {\rm Supp}(f)$. We
have that
 $$(f^{2}\diver
X)(x)=f^{2}(x)[\displaystyle-\frac{v''}{v}(t)+\frac{v'^{2}}{v^{2}}(t)-\frac{v'}{v}(t)\frac{\sqrt{g(t,\xi)}\,'}{\sqrt{g(t,\xi)}}].
$$ Integrating over $B_{N}(p,r)\cap {\rm Supp}(f)$ we have
 \begin{eqnarray}
\int_{B_{N}(p,r)\cap {\rm Supp}(f)}\vert f^{2}\vert \vert\diver
X\vert& \leq & \int_{\xi}\int_{0}^{t
(\xi)}(\displaystyle\vert\frac{v''}{v}\vert+\vert\frac{v'^{2}}{v^{2}}\vert)\vert
f^{2} \vert \sqrt{g(t,\xi)}dtd\xi \nonumber\\
&& \nonumber \\
& &+ \int_{\xi}\int_{0}^{t
(\xi)}\displaystyle\vert\frac{v'}{v}\vert \vert
f^{2}\vert\sqrt{g(t,\xi)} \,'dtd\xi<\infty \nonumber
\end{eqnarray}
that the second term on the right hand side of (\ref{eqDiv(fX)})
is also finite. Where $t(\xi)$ is the largest $t<\min\{d(\xi),\,r
\}$ such that $\exp_{p}(t\,\xi)\in  {\rm Supp}(f)$.  This shows
that $\diver (f^{2}X ) \in L^{1}(B_{N}(p,r))$. Showing the item
(ii)   is a trivial task.

\subsection{Proof of Theorem \protect{\ref{submanifold}}}
 \noin {\bf Theorem \ref{submanifold}} {\em \,Let $N$ be a  Riemannian $n$-manifold  with radial sectional
curvature $K(x)(\partial t, v)\leq c$, for all $x\in
B_{N}(p,r)\setminus {\rm Cut}(p)$, and all $v\perp \partial t $
with $\vert v\vert \leq 1$. Let $M\subset N$ be an $m$-dimensional
minimal submanifold and $\Omega \subset M\cap B_{N}(p,r)$ be a
connected component. Suppose that the $(m-1)$-Hausdorff measure
${\cal H}^{m-1}(\Omega \cap {\rm Cut}_{N}(p))=0$. If $c>0 $,
suppose in addition that $r< \pi /2\sqrt c$. Then
\begin{equation}\label{eqSubm2}\lambda^{\ast}(\Omega) \geq
\lambda_{1}(B_{\mathbb{N}^{m}(c)}(r)),
\end{equation}where $B_{\mathbb{N}^{m}(c)}(r)$  is the geodesic
ball with radius $r$ in the simply connected space form
$\mathbb{N}^{n}(c)$ of constant sectional curvature $c$. If
$\Omega $ is bounded then equality in (\ref{eqSubm1}) holds iff
$\Omega=B _{\mathbb{N}^{m}(c)}(r)$ and $M=\mathbb{N}^{m}(c)$}.

\vspace{.5cm}

 \noin{\bf Proof:}   Let
$v:B_{\mathbb{N}^{m}(c)}(r)\to \mathbb{R}$ be a positive first
Dirichlet eigenfunction of $B_{\mathbb{N}^{m}(c)}(r)$. It is known
that $v$
 is radial  with $v'(t)\leq 0$ and $v'(t)=0$ iff
  $t=0$. We can
normalize $v$ such that $v(0)=1$. The differential equation
$\triangle_{\mathbb{N}^{m}(c)}v(t)+\lambda_{1}(B_{\mathbb{N}^{m}(c)}(r))v(t)=0$
is expressed in geodesic coordinates by
\begin{equation}v''(t) +
(m-1)\frac{C_{c}(t)}{S_{c}(t)}v'(t)+\lambda_{1}(B_{\mathbb{N}^{m}(c)}(r))v(t)=0,
\,\, \forall \,\,t\in[0,r]. \label{eqSubm3}
\end{equation}Recall that for each $\xi \in T_{p}N$, $\vert \xi \vert=1$,  $d(\xi)>0$ is the
largest real  number (possibly $\infty$) such that
  geodesic
  $\gamma_{\xi}(t)=\exp_{p}(t\,\xi)$ minimizes the distance from
  $\gamma_{\xi}(0)=p$ to $\gamma_{\xi}(t)$ for all $t\in [0,d(\xi)]$.
   We have that $B_{N}(p,r)\setminus {\rm Cut}(p)=\exp_{p}(\{t\,\xi\in T_{p}N:
0\le t<\min\{r, d(\xi)\}, \,\vert
   \xi\vert=1 \})$. Define  $u:B_{N}(p,r)\to\mathbb{R}$ by
   $u(\exp_{p}(t\xi))=v(t)$ if $t<\min\{r, d(\xi)\}$ and
   $u(r\xi)=u(d(\xi)\xi)=0$.
 Let
$\Omega \subset M\cap B_{N}(p,r)$ be a connected component and
$\psi:\Omega \rightarrow \mathbb{R}$ defined by $\psi=u\circ
\varphi$, where $\varphi$ is the minimal immersion $\varphi : M
\subset N$. The vector field $X=-\grad \log \psi$  identified with
$d\varphi (X)$  is not smooth at $G= \Omega \cap (\{p\}\cup
 {\rm Cut}(p))$. By hypothesis ${\cal H}^{m-1}(G)=0$
and it can be shown that the vector field $X\in
C^{1}(\Omega\setminus G)\cap L^{\infty}(\Omega)$ and $\diver X\in
L^{1}(\Omega)$  thus $X\in {\cal W}^{1,1}$ and satisfies
(\ref{eqDivergence2}) and by Theorem (\ref{teoremaPrincipal}) we
have that
$$\lambda^{\ast}(\Omega)\geq \inf_{\Omega\setminus G} [\Diver X
-\vert X\vert^{2}]=\inf_{\Omega\setminus G} [\diver X -\vert
X\vert^{2}]=\inf_{\Omega\setminus G}[-\triangle \psi/\psi].$$
Where $\triangle \psi $ is given by the following formula, (see
\cite{kn:bessa-montenegro}, \cite{kn:cheung-leung},
\cite{kn:jorge-koutrofiotis}),
\begin{eqnarray}
\triangle \,\psi (x)    & = & \sum_{i=1}^{m}\Hess\,u(\varphi (x))
\,(e_{i},e_{i})+ \langle \grad u\,,\,
                            \stackrel{\rightarrow}{H}\rangle\label{eqlaplacianpsi}\\
 & =&\sum_{i=1}^{m}\Hess\,u(\varphi (x))
\,(e_{i},e_{i})\nonumber
\end{eqnarray}where $\varphi(x)=\exp_{p}(t\xi) $, $\stackrel{\rightarrow}{H}=0$
 is the mean curvature vector of $\Omega$ at $\varphi (x)$
and $\{ e_{1},\ldots e_{m}\}$ is
an orthonormal basis
 for $T_{\varphi (x)}\,\Omega$. Choose this basis such that $e_{2},\ldots e_{m}$ are
 tangent to the distance sphere $\partial B_{N}(p,t)\subset N$ and
$e_{1}= \cos (\beta(x) )\,\partial/\partial t + \sin (\beta(x)
)\,\partial/\partial \theta $, where $ \partial/\partial \theta
\in [[e_{2}, \ldots e_{m}]]$, $\vert \partial/\partial
\theta\vert=1$. From (\ref{eqlaplacianpsi}) we have for $\varphi
(x)\in \Omega\setminus G$ that

\begin{eqnarray}\triangle \,\psi (x)& =& \sum_{i=1}^{m} \Hess u (\varphi
(x))(e_{i},e_{i})\nonumber \\
&& \nonumber \\
& = & v''(t)(1-\sin^{2}\beta(x) ) \nonumber \\
&&\label{eqlaplacianpsi2} \\
& +& v'(t)\,\sin^{2}\beta(x)\,{\rm
Hess}(t)(\partial /\partial\theta ,\partial/\partial \theta ) \nonumber \\
&&\nonumber \\
 & +& v'(t)\,\sum_{i=2}^{m}{\rm Hess}(t)(e_{i},e_{i})  \nonumber
\end{eqnarray}where $t={\rm dist}_{N}(p,x)$.
Adding and subtracting $(C_{c}/S_{c})(t)\,v'(t)\sin^{2}\beta(x) $
and $ (m-1)(C_{c}/S_{c})(t)\,v'(t)$ in (\ref{eqlaplacianpsi2}) we
have
\begin{eqnarray}\label{eqlaplacianpsi3}
\triangle \,\psi (x)&=&v''(t)+ \,(m-1)\frac{C_{c}}{S_{c}}(t) \,v'(t)\nonumber \\
&&\nonumber \\
& +& \left({\rm Hess}(t)(\partial /\partial\theta
,\partial/\partial \theta )-\frac{C_{c}(t)}{S_{c}(t)}\right)
\,v'(t)\sin^{2}\beta(x) \nonumber \\
&& \\
&+& \sum_{i=2}^{m}[{\rm Hess
}(t)(e_{i},e_{i})-\frac{C_{c}}{S_{c}}(t)]\,v'(t)\nonumber \\
&&\nonumber \\
& +&  \nonumber\left(
\frac{C_{c}(t)}{S_{c}(t)}v'(t)-v''(t)\right)\sin^{2}\beta(x)
\end{eqnarray}
From (\ref{eqSubm3}) and (\ref{eqlaplacianpsi3}) we have that
\begin{eqnarray}\displaystyle-\frac{\triangle \,\psi}{\psi} (x)& =&
\label{eqlaplacianpsi4} \lambda_{1}(B_{\mathbb{N}^{m}(c)}(r))\nonumber \\
& -& \left({\rm Hess}(t)(\partial /\partial\theta
,\partial/\partial \theta )-\frac{C_{c}(t)}{S_{c}(t)}\right)
\,\frac{v'(t)}{v(t)}\sin^{2}\beta(x) \nonumber \\
&& \\
&-& \sum_{i=2}^{m}[{\rm Hess
}(t)(e_{i},e_{i})-\frac{C_{c}}{S_{c}}(t)]\,\frac{v'(t)}{v(t)}\nonumber \\
&&\nonumber \\
& -&\frac{1}{v(t)} \left(
\frac{C_{c}(t)}{S_{c}(t)}v'(t)-v''(t)\right)\sin^{2}\beta(x)
\nonumber
\end{eqnarray}Since the radial curvature  $K(x)(\partial t, v)\leq c$ for all $x\in
B_{N}(p,r)\setminus {\rm Cut}(p)$ and all $v\perp \partial t $
with $\vert v\vert \leq 1$ then by the Hessian Comparison Theorem
(see \cite{kn:schoen-yau}) we have that $\Hess (t(x))(v,v) \geq
(C_{c}/S_{c})(t)$ for all $v\perp \partial t$, $t(x)=t$,
$x=\exp_{p}(t\xi)$. But $v'(t)\leq 0$ then we have that the second
and third terms of (\ref{eqlaplacianpsi4}) are non-negative. If
the fourth term of (\ref{eqlaplacianpsi4}) is non-negative then we
would have that
$$-\frac{\triangle \,\psi }{\psi} (x)  \geq
\lambda_{1}(B_{\mathbb{N}^{m}(c)}(r)) $$ \noin By Theorem
(\ref{teoremaPrincipal}) we have that
\begin{equation}\label{eqlaplacianpsi5}\lambda^{\ast}(\Omega
)\geq \inf (-\frac{\triangle \,\psi }{\psi} )\geq
\lambda_{1}(B_{\mathbb{N}^{m}(c)}(r)). \end{equation} This proves
(\ref{eqSubm2}). We can see that $\displaystyle -\left(
\frac{C_{c}(t)}{S_{c}(t)}\frac{v'(t)}{v(t)}-\frac{v''(t)}{v(t)}\right)\sin^{2}\beta(x)\geq
0$ is equivalent to
\begin{equation}\label{eqSubm4}m
\frac{C_{c}(t)}{S_{c}(t)}v'(t)+\lambda_{1}(B_{\mathbb{N}^{m}(c)}(r))v(t)
< 0, \,\, t\in (0,r).
\end{equation}
To prove (\ref{eqSubm4}) we will assume without loss of generality
that $c=-1,0,1$. Let us consider first the case $c =0$ that
presents the idea of the proof.  The other two remaining cases
($c=-1 $ and $c=1$) we are going to treat (quickly) with the same
idea.
 When $c=0$ the inequality (\ref{eqSubm4}) becomes
\begin{equation} \label{eqSubm5}\frac{mv'(t)}{t}+ \lambda_{1}v(t) <
0, \,\, t\in (0,r),
\end{equation}where $\lambda_{1}:=\lambda_{1}(B_{\mathbb{N}^{m}(c)}(r))$.
Let $\mu (t):=\exp\{-\displaystyle\frac{\lambda_{1}t^{2}}{2m}\}$.
The functions $v$ and $\mu$ satisfy the following identities,
\begin{equation}\label{eqSubm6}\begin{array}{lcl} (t^{m-1}v'(t))' +
\lambda_{1}t^{m-1}v(t)& =& 0\\
&& \\
(t^{m-1}\mu'(t))' + \lambda_{1}t^{m-1}(1-
\displaystyle\frac{\lambda_{1}\,t^{2}}{m^{2}})\mu(t)& =& 0
\end{array}
\end{equation}In (\ref{eqSubm6}) we multiply the first identity by
$\mu$ and the second by $-v$ adding them and integrating from $0$
to $t$ the resulting identity we obtain, $$t^{m-1}\,v'(t)\,\mu(t)
-t^{m-1}\,v(t)\,\mu'(t)=-\frac{\lambda_{1}^{2}}{m^{2}}\int_{0}^{t}\mu(t)\,v(t)
<0,\,\,\ \forall t\in (0,r).
$$
Then $\mu(t)v'(t) < \mu'(t)v(t) $ and this proves (\ref{eqSubm5}).

\vspace{.1cm}

 \noin  Assume that now that $c=-1$. The inequality
(\ref{eqSubm4}) becomes
\begin{equation}\label{eqSubm7}m\displaystyle\frac{C_{-1} (t)}{S_{-1}(t)}v'(t) +
\lambda_{1}v(t)<0
\end{equation}
Set $\mu(t):=C_{-1} (t)^{-\lambda_{1}/m}$. The functions $v$ and
$\mu$ satisfy the the following identities
\begin{equation}\label{eqSubm8}\displaystyle\begin{array}{lll}(S_{-1}^{m-1}v')' +
\lambda_{1}S_{-1}^{m-1}v& =& 0 \\
&& \\
 (S_{-1} ^{m-1}\mu')' +
 \lambda_{1}S_{-1}^{m-1}\left(\displaystyle \frac{m-1}{m}+\displaystyle \frac{1}{m
C_{-1}^{2} } -
\displaystyle\frac{\lambda_{1}}{m^{2}}\displaystyle\frac{S_{-1}^{2}}{C_{-1}^{2}}\right)\mu&=&0
\end{array}
\end{equation}In (\ref{eqSubm8}) we multiply the first identity by
$\mu$ and the second by $-v$ adding them and integrating from $0$
to $t$ the resulting identity we obtain
$$S_{-1}^{m-1}\left(v'\mu-\mu'v\right)(t)+
\int_{0}^{t}\lambda_{1}S_{-1}^{m-1}\left[\frac{1}{m}-\frac{1}{mC_{-1}^{2}}+
\frac{\lambda_{1}}{m^{2}}\frac{S_{-1}^{2}}{C_{-1}^{2}}\right]\mu
v=0$$ The  term $
S_{-1}^{m-1}\left[\displaystyle\frac{1}{m}-\frac{1}{mC_{-1}^{2}}+
\frac{\lambda_{1}}{m^{2}}\frac{S_{-1}^{2}}{C_{-1}^{2}}\right]\mu v
$  is positive (one can easily check) therefore we have that
$(v'\mu-\mu'v)(t)<0$ for all $t\in (0,r)$. This  proves
(\ref{eqSubm7}).

\vspace{.1cm}

 \noin  For $c=1$ the inequality (\ref{eqSubm4}) becomes
 the following inequality
\begin{equation}\label{eqSubm9}m\frac{C_{1}}{S_{1}}v'(t)+\lambda_{1}v(t)<0,\,\,\,0<t<\pi/2
\end{equation}Set
$  \mu(t):=C_{1} (t)^{-\lambda_{1}/m},\,\,0<t<\pi/2$. The
functions $v$ and $\mu$ satisfy the the following identities
\begin{equation}\label{eqSubm10}\begin{array}{lll}(S_{1}^{m-1}v')' +
\lambda_{1}S_{1}^{m-1}v&=& 0\\
&& \\
 (S_{1} ^{m-1}\mu')' -
 \lambda_{1}S_{1}^{m-1}\left(\displaystyle \frac{m-1}{m}+ \displaystyle\frac{1}{m
C_{1}^{2} } +
\displaystyle\frac{\lambda_{1}}{m^{2}}\displaystyle\frac{S_{1}^{2}}{C_{1}^{2}}\right)\mu&=&0
\end{array}
\end{equation}In (\ref{eqSubm10}) we multiply the first identity by
$\mu$ and the second by $-v$ adding them and integrating from $0$
to $t$ the resulting identity we
obtain$$S_{1}^{m-1}\left(v'\mu-\mu'v\right)(t)+
\int_{0}^{t}\lambda_{1}S_{1}^{m-1}\left[2-\frac{1}{m}+\frac{1}{mC_{1}^{2}}+
\frac{\lambda_{1}}{m^{2}}\frac{S_{1}^{2}}{C_{1}^{2}}\right]\mu
v=0$$ The term $\displaystyle
S_{1}^{m-1}\left[2-\frac{1}{m}+\frac{1}{mC_{1}^{2}}+
\frac{\lambda_{1}}{m^{2}}\frac{S_{1}^{2}}{C_{1}^{2}}\right]\mu v$
is positive therefore we have that $(v'\mu-\mu'v)(t)<0$ for all
$t\in (0,r)$. This  proves (\ref{eqSubm9}) and thus the fourth
term in (\ref{eqlaplacianpsi4}) is non-negative.

To finishes the proof of the Theorem (\ref{submanifold}) we need
to consider the equality case in (\ref{eqSubm2}) when $\Omega $ is
bounded. From the spectral theory it is known that for a given
bounded domain in a Riemannian manifold there is $u\in
C^{\infty}(\Omega)\cap H_{0}^{1}(\Omega)$, positive in $\Omega$
satisfying $\triangle u+\lambda_{1}(\Omega)u=0$, where
$\lambda_{1}(\Omega)=\lambda^{\ast}(\Omega)$. This function is
also called an {\em eigenfunction}. Observe that $u\vert
\partial \Omega =0$ only a.e. thus $u$ is not considered  a solution
for the Dirichlet eigenvalue problem. The proof of existence of
$u$ is the same proof of existence of Dirichlet eigenvalues for a
smooth domain, since  it (the proof) does need smoothness of the
boundary but the boundedness of the domain.  With an approximation
argument Barta's Theorem can be extend to arbitrary bounded open
sets.
\begin{proposition}\label{3.1}Let $\Omega$ be a bounded domain in a smooth Riemannian manifold. Let
$v\in C^{2}(\Omega)\cap C^{0}(\overline{\Omega})$, $v>0$ in
$\Omega$ and $v\vert
\partial \Omega =0$. Then \begin{equation}\lambda^{\ast}(\Omega)\geq
\inf_{\Omega}(-\frac{\triangle
v}{v}).\label{eqProp1}\end{equation}Moreover,
$\lambda^{\ast}(\Omega)=
\inf_{\Omega}(-\displaystyle\frac{\triangle v}{v})$ if and only if
$v=u$,  where $u$ is a positive eigenfunction of $\Omega$, i.e.
$\triangle u+\lambda^{\ast}(\Omega)u=0$.
\end{proposition}
\noin {\bf Proof:}  Let $\epsilon_{i}\to 0$ be a sequence of
positive  regular values of $v$ and let
$\Omega_{\epsilon_{i}}^{v}=\{x\in \Omega;\, v(x)>\epsilon_{i}\}$.
Applying Barta's Theorem we have that
\begin{equation}\label{eqProp2}\lambda^{\ast}(\Omega_{\epsilon_{i}}^{v})=
\lambda_{1}(\Omega_{\epsilon_{i}}^{v})\geq
\inf_{\Omega_{\epsilon_{i}}^{v}}(-\displaystyle\frac{\triangle
v}{v})\geq \inf_{\Omega}(-\displaystyle\frac{\triangle
v}{v})\end{equation}But $\lim_{\epsilon_{i}\to
0}\lambda^{\ast}(\Omega_{\epsilon_{i}}^{v})=\lambda^{\ast}(\Omega)$,
see in \cite{kn:chavel}, page 23.  Let $u\in C^{\infty}\cap
H_{0}^{1}(\Omega)$ be a positive eigenfunction of $\Omega$. Then
\begin{equation}\label{eqProp3}\lambda^{\ast}(\Omega)= \displaystyle
-\frac{\triangle u}{u} =-\frac{\triangle v}{v}+\frac{u\triangle
v-v\triangle u}{uv}
\end{equation}
\noin {\bf Claim:} $\displaystyle\int_{\Omega}(u\triangle
v-v\triangle u)=0$.\\

\noin {\bf Proof:} Let $\epsilon_{i}\to 0$ be a sequence of
positive  regular values of $u$ and let
$\Omega_{\epsilon_{i}}^{u}=\{x\in \Omega;\, u(x)>\epsilon_{i}\}$.
Set
\begin{equation}\label{eqProp4}u_{\epsilon_{i}}=\left\{\begin{array}{lll}
u-\epsilon_{i} & {\rm on} & \Omega_{\epsilon_{i}}^{u} \\
&& \\
0 &{\rm on} & \Omega
\setminus\Omega_{\epsilon_{i}}^{u}\end{array}\right.
\end{equation}One can show  that $u_{\epsilon_{i}}\to u$ in
$H^{1}(\Omega)$ using the Lebesgue Convergence Theorem. Therefore
\begin{eqnarray}
\displaystyle\int_{\Omega_{\epsilon_{i}}^{u}}(u_{\epsilon_{i}}\triangle
v-v\triangle u_{\epsilon_{i}})
&=&-\displaystyle\int_{\Omega_{\epsilon_{i}}^{u}}\langle \grad
u_{\epsilon_{i}}, \grad v\rangle +
\lambda^{\ast}(\Omega)\int_{\partial\Omega_{\epsilon_{i}}^{u}}u
v\nonumber \\&& \\
 & \to &\displaystyle-\int_{\Omega}\langle \grad u, \grad v\rangle +
\lambda^{\ast}(\Omega) \int_{\partial\Omega}u v=0\nonumber
\end{eqnarray}Since  $v\in H_{0}^{1}(\Omega)$ and $u$ is  a
weak solution of $\triangle u+\lambda^{\ast}(\Omega)u=0$. On the
other hand   $
\displaystyle\int_{\Omega_{\epsilon_{i}}^{u}}(u_{\epsilon_{i}}\triangle
v-v\triangle u_{\epsilon_{i}})\to
\displaystyle\int_{\Omega}(u\triangle v-v\triangle u) $.

 \noin If  $\Omega$ is bounded  we have that $\partial \varphi (\Omega)\subset \partial
 B_{N}(p,r)$. This implies that the function $\psi= u\circ
 \varphi\in C^{2}(\Omega)\cap C^{0}(\overline{\Omega})$ is such that
 $\psi \,\vert \partial \Omega =0$. Suppose that
 $ \lambda^{\ast}(\Omega)=\lambda_{1}(B_{\mathbb{N}^{m}(c)}(r))$.
 Then  by Proposition (\ref{3.1}),  $\psi:\Omega\to
\mathbb{R}$ is an eigenfunction of $\Omega$ and we have that $
\lambda^{\ast}(\Omega)=-\triangle \psi/\psi$. From
(\ref{eqlaplacianpsi4}) we have that
\begin{equation}\begin{array}{lll}  \left({\rm Hess}(t)(\partial /\partial\theta
,\partial/\partial \theta
)-\displaystyle\frac{C_{c}(t)}{S_{c}(t)}\right)
\,\displaystyle\frac{v'(t)}{v(t)}\sin^{2}\beta(x)&=&0  \\
&& \\
\displaystyle\sum_{i=2}^{m}[{\rm Hess
}(t)(e_{i},e_{i})-\displaystyle\frac{C_{c}}{S_{c}}(t)]\,\displaystyle\frac{v'(t)}{v(t)}&=&0 \\
&&\nonumber \\
\displaystyle\frac{1}{v(t)} \left(
\frac{C_{c}(t)}{S_{c}(t)}v'(t)-v''(t)\right)\sin^{2}\beta(x) &=&0,
\end{array}
\end{equation}for all $t$ such that $\varphi (x)=\exp_{p}(t\xi)\in \Omega$.
This implies $\sin^{2}\beta (x)=0$ for all $x\in \Omega$ and we
have that $e_{1}(\varphi (x))=\partial /\partial t $. Integrating
the vector field $\partial/\partial t$ we have a minimal geodesic
(in $N\cap \varphi ( \Omega)$) joining $\varphi (x)$ to the center
$p$. This imply that $\Omega$ is the geodesic ball in $M$ centered
at $\varphi^{-1}(p)$ with radius $r$ i.e. $\Omega
=B_{M}(\varphi^{-1}(p),r)$. Since $\psi$ is an eigenfunction with
the same eigenvalue $\lambda_{1}(B_{\mathbb{N}^{m}(c)}(r))$ we
have that
\begin{equation}\label{eq29}\triangle_{M}\,v (t)=\triangle_{\mathbb{N}^{m}(c)}\,v(t),\,\,t=dist_{N}(p,\varphi
(q)), \forall q\in \Omega .
\end{equation}Rewriting this identity (\ref{eq29}) in geodesic  coordinates we have that
$$\frac{\sqrt{g(t,\xi)}'}{\sqrt{g(t,\xi)}} (t,\theta)v'(t) + v''(t)=
(m-1)\frac{C_{c}(t)}{S_{c}(t)}v'(t) + v''(t) $$ \noin This imply
that by Bishop Theorem $\Omega= B_{M}(\varphi^{-1}(p)$ and
$B_{\mathbb{N}^{m}(c)}(r)$ are isometric. By analytic continuation
$M=\mathbb{N}^m(c)$.

\subsection{Proof of Corollary \ref{nadirashvili}}Let $M\subset
B_{\mathbb{R}^{3}}(r)$ be a
 complete bounded minimal surface in $\mathbb{R}^{3}$. Then
 $$\lambda^{\ast}(M)\geq
 \lambda_{1}(\mathbb{D}(r))=c/r^{2}.
$$ \noin
 Where $c>0$ is an absolute constant. The proof of this result
 follows directly from Theorem (\ref{submanifold}). The theorem
 says that the fundamental tone $\lambda^{\ast}(\Omega)\geq
 \lambda_{1}(B_{\mathbb{R}^{3}}(r))$ for any connected component
 of $M\cap B_{\mathbb{R}^{3}}(r)$. In particular,  for $\Omega= M$
 we have that $\lambda^{\ast}(M)\geq
 \lambda_{1}(B_{\mathbb{R}^{3}}(r))=c(2)/r^{2}$, where $c(2)$ is
 the first zero of the Bessel function $J_{0}$, see
 \cite{kn:chavel}, page 46.

\subsection{Proof of Corollary \protect{\ref{stability}}}

Let $\varphi :M \hookrightarrow N^{n+1}$ be a complete orientated
minimal hypersurface and $A(X) =- \nabla_{X}\, \eta$ its  second
fundamental form, where $\eta$ is globally defined unit vector
field normal to $\varphi (M)$. A normal domain  $D\subset M$ is
said to be stable if  the first Dirichlet eigenvalue
$\lambda_{1}^{L}=\inf\{-\smallint_{D} u\,Lu/\smallint_{D}
u^2,\,u\in C_{0}^{\infty}(D)\}$ of the operator $L=\triangle +
Ric(\eta)+ \Vert A\Vert^{2}$ is positive.  On the other hand we
have that $$\begin{array}{lll}-\smallint_{D}u\,Lu & = &
\smallint_{D} \left[\,\vert \grad u
\vert^{2}-\left( Ric (\eta)+\Vert A\Vert^{2}\right)u^{2}\right]\\
&& \\
 & \geq &\smallint_{D}\lambda_{1}^{\triangle}(D)- \left( Ric
(\eta)+\Vert A\Vert^{2}\right)u^{2}.
\end{array}
$$ \noin Therefore if $\lambda_{1}^{\triangle}(D)\geq \sup_{x\in D}\{ Ric
(\eta)+\Vert A\Vert^{2}(x)\}$ then $D$ is stable. In
\cite{kn:bessa-montenegro} we give estimates for
$\lambda_{1}^{\triangle}(D)$ in submanifolds with locally bounded
mean curvature, in particular minimal hypersurfaces. With those
estimates we have obvious statements for stability theorems.  If
$D=B_{M}(p,r)$ is a ball in $M^{n}\hookrightarrow \mathbb{R}^{n+1}
$ obviously that $\varphi (D)\subset B_{n+1}(\varphi (p),r)$. And
we have that $\lambda_{1}(D)\geq
\lambda^{\ast}(\varphi^{-1}(B_{n+1}(\varphi (p),r))\geq
\lambda_{1}(B_{n}(0,r))$. This proves the Corollary
(\ref{stability}).

\subsection{Quasilinear elliptic equations}
\noin In this section we want to apply Barta's Theorem to study
the existence of solutions to certain quasi-linear elliptic
equations. Let $M$ be a bounded Riemannian manifold with smooth
non-empty boundary $\partial M$ and $f\in C^{2}(M)$, $f>0$. If we
set $u=-\log f$\footnote{ This transformation $u=-\log f$ we
learned from Kazdan \& Kramer \cite{kn:kazdan-kramer} but  it also
appears
 in \cite{kn:fischer-colbrie-schoen}}, $f>0$ then the problem $ \triangle\,f + F f=0$
becomes $\triangle u - \vert \grad u\vert^{2}= F $. Hence Barta's
Theorem (\ref{barta}) can be translated as

\vspace{.1cm}

 \noin {\bf Theorem \ref{Elliptic1} } \, {\em Let $M$
be a bounded Riemannian manifold with smooth boundary and
 $F\in C^{0}(\overline{M})$.
Consider this problem,
\begin{equation}\label{eqElliptic3}\left\{\begin{array}{rccrl} \triangle u - \vert
\grad u\vert^{2} &=& F & in&
M\\
u&=& +\infty & on& \partial M .
\end{array}\right.
\end{equation}If (\ref{eqElliptic3})  has a smooth solution then $\inf_{M} F\leq
\lambda_{1}(M) \leq \sup_{M} F $. If either $\inf_{M}
F=\lambda_{1}(M)$ or $\lambda_{1}(M)= \sup_{M} F$ then
$F=\lambda_{1}(M)$.   On the other hand if
 $F=\lambda$ is a constant
 the problem (\ref{eqElliptic3}) has solution if and only if
$\lambda=\lambda_{1}(M)$.}

\vspace{.3cm} Likewise Theorem (\ref{teoremaPrincipal}) can be
translated into language of quasi-linear elliptic equations.
\vspace{.2cm}

 \noin {\bf Theorem \ref{Elliptic2}}\,{\em  Let $M$
be a bounded Riemannian manifold with smooth boundary and
 $F\in C^{0}(\overline{M})$ and $\psi \in C^{0}(\partial
M)$. Consider the problem
\begin{equation}\label{eqElliptic4}\left\{\begin{array}{rclll} \triangle u
- \vert \grad u\vert^{2} &=& F & in&
M\\
u&=& \psi & on& \partial M.
\end{array}\right.
\end{equation}then
~

\begin{enumerate}
\item[a)] If $\sup_{M}F < \lambda_{1}(M)$ then (\ref{eqElliptic4})
has solution. \item[b)] If (\ref{eqElliptic4}) has
 solution then $\inf_{M}F< \lambda_{1}(M)$.
\end{enumerate}}

\noin {\bf Proof:} The operator $L=-\triangle -F $ is compact thus
its spectrum is a sequence of eigenvalues $\lambda_{1}^{L}<
\lambda_{2}^{L}\leq \lambda_{3}^{L}\leq \cdots \nearrow +\
\infty$.  Suppose that $\sup_{M}F<\lambda_{1}(M)$ then
 we
have that \begin{eqnarray}\lambda_{1}^{L}& = & \inf\,\{ \int_{M}
\vert \grad u\vert^{2} -F\,u^{2};\, u\in
H_{0}^{1}(M),\,\smallint_{M}u^{2}=1\}\nonumber \\
 & \geq & \lambda_{1}(M) - \sup_{M}F>0.\nonumber
\end{eqnarray} Therefore $L$ is invertible. On the other hand,
$u=-\log f $ is a solution of (\ref{eqElliptic4}) if and only if
$f$ is solution of
\begin{equation}\label{eqElliptic5}\left\{\begin{array}{lcl}
L f & = & 0 \,\,\,\, in \,\, M \\
f &= & e^{-\psi}\,\, on \, \, \partial M
\end{array}\right.
\end{equation} Consider the harmonic extension $v$ of $e^{-\psi}$ on $\partial
M$, ($\triangle \,v=0$ in $M $ and $v=e^{-\psi}$ on $\partial M$),
and $h= f-v$.  We have that $L \,h = F\, v$ in $M $ and $h=0 $ on
$\partial M$. Then $h = L^{-1}(F\,v)$ and $f= -v + L^{-1} (F \,v)$
is solution of  (\ref{eqElliptic5}).

\vspace{.1cm}

 \noin  Item b. Suppose that the equation
(\ref{eqElliptic4}) has a smooth solution $u$. Let $\varphi$ be
the first eigenfunction of the operator $-\triangle $,
($\triangle\, \varphi + \lambda_{1}(M ) \varphi =0$ in $M $ and
$\varphi\vert\partial M =0$), and $f$ a solution of
(\ref{eqElliptic5}). By Green we have that
$$\int_{ M }f \triangle \varphi - \varphi \triangle f=
\int_{\partial M }e^{-\psi}\frac{\partial \varphi}{\partial \eta}-
\varphi \frac{\partial f}{\partial \eta}= \int_{\partial
M}e^{-\psi}\frac{\partial \varphi}{\partial \eta}\;.
$$
Thus $$\int_{\partial M}e^{-\psi}\frac{\partial \varphi}{\partial
\eta} = \int_{M}(F-\lambda_1)\varphi f $$ if $v$ is the harmonic
extension $v$ of $e^{-\psi}$ then we have
$$\int_{M}-\lambda_1f\varphi =\int_{ M }v \triangle
\varphi - \varphi \triangle v= \int_{\partial M }v\frac{\partial
\varphi}{\partial \eta}- \varphi \frac{\partial v}{\partial
\eta}=\int_{\partial M}e^{-\psi}\frac{\partial \varphi}{\partial
\eta}<0$$
 since
$\inf_{M}v=\inf_{\partial M}v=\inf_{\partial M}e^{-\psi}>0$.
Therefore $\smallint_{M}(F-\lambda_1)\varphi f <0$ and $\inf
F<\lambda_1(M)$.


\begin{thebibliography}{abcd}
\bibitem{kn:barta} Barta, J.: {\em Sur la vibration fundamentale
d'une membrane.} C. R. Acad. Sci. 204, 1937, 472-473.
\bibitem{kn:berger-gauduchon-mazet} Berger, M., Gauduchon, P. and Mazet,
E.:  {\em Le Spectre d'une Vari\'{e}t\'{e} Riemannienes}. Lect.
Notes Math. 194, 1974, Springer-Verlag.
\bibitem{kn:bessa-montenegro} Bessa, G. P., Montenegro, J. F.:
{\em Eigenvalue estimates for submanifolds  with locally bounded
mean  curvature.} Ann.  Global Anal. and Geom.  24, 2003, 279-290.
\bibitem{kn:bishop-crittenden} Bishop, R., Crittenden, R.: {\em
Geometry of Manifolds}. Academic Press, New York, 1964.
\bibitem{kn:chavel} Chavel, I.: {\em Eigenvalues in Riemannian
Geometry.}  Pure and Applied Mathematics, 1984, Academic Press,
INC.
\bibitem{kn:cheng1} Cheng, S. Y., {\em Eigenfunctions and eigenvalues of the Laplacian.}
 Am. Math. Soc Proc. Symp. Pure Math. 27, part II, 1975, 185-193.
 \bibitem{kn:cheung-leung} Cheung, Leung-Fu and  Leung, Pui-Fai, {\em Eigenvalue estimates for submanifolds with bounded
 mean curvature in the hyperbolic space.} Math. Z.  236, 2001, 525-530.
\bibitem{kn:cheng-li-yau} Cheng, S. Y., Li, P. and Yau, S. T., {\em
Heat equations on minimal submanifolds and their applications.}
Amer. J. Math. 106, 1984, 1033-1065.
\bibitem{kn:fischer-colbrie-schoen} Fisher-Colbrie, D. and Schoen, R,.
{\em The structure of complete stable minimal surfaces in
$3$-manifolds of non negative scalar curvature.}  Comm.  Pure and
Appl. Math. 33, 1980, 199-211.
\bibitem{kn:G} Grigor'yan, A.: {\em Analytic and geometric
 background of recurrence and non-explosion of the brownian motion
 on Riemannian manifolds}. Bull. Amer. Math. Soc. 36, 2, 1999,
 135-249.
\bibitem{kn:kazdan-kramer} Kazdan, J. and Kramer, R., {\em
Invariant Criteria for existence of solutions to second-order
quasilinear elliptic equations.} Comm. Pure Appl. Math. 31, 5,
1978, 619-645.
\bibitem{kn:kazdan-warner}Kazdan, J. and Warner, F. W., {\em
Remarks on some quasilinear elliptic equations.} Comm. Pure Appl.
Math. 28, 1975, 567-597.
\bibitem{kn:jorge-koutrofiotis} Jorge, L. and  Koutrofiotis, D., {\em An estimate for the curvature of
bounded submanifolds.} Amer. J. Math., 103, 4, 1980, 711-725.
\bibitem{kn:mckean}McKean, H. P.:  {\em An upper bound for the spectrum of $\triangle $ on a manifold of
 negative curvature.} J. Differ. Geom. 4, 1970, 359-366.
\bibitem{kn:matilla} Matilla, P., {\em Geometry of Sets and
Measures in Euclidean Spaces.} Cambridge University Press, 1995.
\bibitem{kn:Nadirashvili} Nadirashvili, N., {\em Hadamard's and
Calabi-Yau's conjectures on negatively curved and minimal
surfaces.} Invent. Math., 126, 1996, 457-465.
\bibitem{kn:schoen} Schoen, R., {\em Stable minimal surfaces in
three manifolds.}  Seminar on Minimal Submanifolds. Annals of Math
Studies. Princenton  University Press.
\bibitem{kn:schoen-yau} Schoen, R. and Yau, S. T.,  {\em Lectures on Differential Geometry.}
Conference Proceedings and Lecture Notes in Geometry and Topology,
{vol. 1}, 1994.
\bibitem{kn:yau} Yau, S. T., {\em Review of Geometry and
Analysis.} Asian J. Math. vol.4, 1, 2000, 235-278.
\bibitem{kn:whitney} Whitney, H., {\em Geometric Integration
Theory.} Princenton Mathematical Series, 1957, Princenton
University Press.

\end{thebibliography}
\end{document}